\documentclass[a4paper,11pt]{amsart}
\usepackage{graphicx}
\usepackage{mathptmx}
\usepackage{amsmath}
\usepackage{amssymb}
\usepackage{enumitem}
\usepackage{xcolor}

\newmuskip\pFqmuskip

\newcommand*\pFq[6][8]{%
  \begingroup 
  \pFqmuskip=#1mu\relax
  \mathcode`=\string"8000
  \begingroup\lccode`\~=`\,
  \lowercase{\endgroup\let~}\pFqcomma
  F^{#2}_{#3}{\left(\genfrac..{0pt}{}{#4}{#5}\bigg|#6\right)}%
  \endgroup
}
\newcommand{\pFqcomma}{\mskip\pFqmuskip}

\newtheorem{theorem}{Theorem}

\newtheorem{corollary}[theorem]{Corollary}

\begin{document}

\title[Poisson degenerate central moments]{Poisson degenerate central moments related to degenerate Dowling and degenerate $r$-Dowling polynomials}

\author{Taekyun Kim$^{1, *}$}
\address{Kwangwoon University, Seoul, 139-701, Republic of Korea}
\email{tkkim@kw.ac.kr}

\author{Dae San Kim$^{2, *}$}
\address{Sogang University, Seoul, 121-742, Republic of Korea}
\email{dslee@daegu.ac.kr}

\author{Hye Kyung  Kim$^{3,*}$}
\address{Department of Mathematics Education, Daegu Catholic University, Gyeongsan 38430, Republic of Korea}
\email{hkkim@cu.ac.kr}

\subjclass[2010]{11B73; 11B83}
\keywords{degenerate $r$-Whitney numbers of the second kind; degenerate $r$-Dowling polynomials; Poisson random variable;  Poisson degenerate central moments; Charlier polynomials}
\thanks{* are corresponding authors}

\begin{abstract}Degenerate Dowling and degenerate $r$-Dowling polynomials were introduced earlier as degenerate versions and further generalizations of Dowling and $r$-Dowling polynomials. The aim of this paper is to show their connections with Poisson degenerate central moments for a Poisson random variable with a certain parameter and with Charlier polynomials.
\end{abstract}

\maketitle

\markboth{\centerline{\scriptsize Poisson degenerate central moments related to degenerate Dowling and degenerate $r$-Dowling polynomials}}
{\centerline{\scriptsize  T. Kim, D. S. Kim and  H. K. Kim}}

\section{introduction and preliminaries}
In recent years, studying various degenerate versions of many special polynomials and numbers, which began with the paper [4]  by Carlitz, received regained interests of some mathematicians and many interesting results were discovered (see\,\,[12-18] and the references therein). They have been explored by employing several different tools such as combinatorial methods, generating functions, $p$-adic analysis, umbral calculus techniques, differential equations, probability theory and analytic number theory.\par
Degenerate Dowling and degenerate $r$-Dowling polynomials were introduced earlier as degenerate versions and further generalizations of Dowling and $r$-Dowling polynomials. The aim of this paper is to show their connections with Poisson degenerate central moments for a Poisson random variable with a certain parameter and with Charlier polynomials. \par
The outline of this paper is as follows. In Section 1, we recall the Stirling numbers of the first and second kinds, Bell polynomials, the degenerate exponential functions, the degenerate Stirling numbers of the first and second kinds, the degenerate Bell polynomials, the Poisson random variable with parameter $\alpha$, and the Charlier polynomials. In addition, we remind the reader of the Whitney numbers of the first and second kinds, Dowling polynomials, the degenerate Whitney numbers of the second kind, the degenerate Dowling polynomials, the degenerate $r$-Whitney numbers of the second kind and the degenerate $r$-Dowling polynomials. Section 2 is the main result of this paper. In the following, assume that $X$ is the Poisson random variable with mean $\frac{\alpha}{m}$. In Theorem 1, we show that  the Poisson degenerate central moment $E[(mX+1)_{n,\lambda}]$ is equal to $D_{m,\lambda}(n,\alpha)$, where $D_{m,\lambda}(n,x)$ is the degenerate Dowling polynomial. In Theorem 2, we express the same Poisson degenerate central moment in terms of the degenerate Bell polynomials. In Theorem 4, we deduce that $E[(mX+r)_{n,\lambda}]$ is equal to $D_{m,\lambda}^{(r)}(n,\alpha)$, where $D_{m,\lambda}^{(r)}(n,x)$ is the degenerate $r$-Dowling polynomial. We express the same in terms of the degenerate Bell polynomials in Corollary 5, and of the degenerate $r$-Whitney numbers of the second kind and the Bell polynomials in Theorem 6. Furthermore, it is represented by the Charlier polynomials and the degenerate Stirling numbers of the second kind in Theorems 10 and 11. In the rest of this section, we recall the facts that are needed throughout this paper.

\par

It is well known that the stirling numbers of the first kind are defined as
\begin{equation}\label{eq01}
\begin{split}
(x)_n = \sum_{k=0}^n S_1(n,k)x^k, \ \ \ (n\geq0), \quad \ \ \text{(see \ \cite{1, 2, 3, 6, 7, 21})},
\end{split}
\end{equation}
where $(x)_0=1, \ (x)_n=x(x-1) \cdots (x-n+1), \ \ (n\geq1).$

The Stirling numbers of the second kind are given by
\begin{equation}\label{eq02}
\begin{split}
x^n=\sum_{k=0}^n S_2(n,k) (x)_k, \ \ \ (n\geq0), \quad \ \ \text{(see \ \cite{18, 19})}.
\end{split}
\end{equation}

From \eqref{eq01} and \eqref{eq02}, we note that

\begin{equation}\label{eq03}
\begin{split}
\frac{1}{k!} \big(\log(1+t)\big)^k = \sum_{n=k}^\infty S_1(n,k)\frac{t^n}{n!},
\end{split}
\end{equation}

and

\begin{equation}\label{eq04}
\begin{split}
\frac{1}{k!}(e^t-1)^k=\sum_{n=k}^\infty S_2(n,k)\frac{t^n}{n!},\quad (k\geq0),\quad \ \text{(see \ \cite{5, 6, 17, 19})}.
\end{split}
\end{equation}

The Bell polynomials are defined by
\begin{equation}\label{eq05}
\begin{split}
e^{x(e^t-1)}=\sum_{n=0}^\infty \phi_n(x)\frac{t^n}{n!},\quad \ \text{(see \ \cite{5, 6, 19})}.
\end{split}
\end{equation}

Thus, by \eqref{eq05}, we get
\begin{equation}\label{eq06}
\begin{split}
\phi_n(x)=\sum_{k=0}^n S_2(n,k)x^k, \quad (n\geq0), \quad \ \text{(see \ \cite{4, 16, 19})}.
\end{split}
\end{equation}

In \cite{12}, the degenerate exponentials are defined by
\begin{equation}\label{eq07}
\begin{split}
e_\lambda^x(t)=(1+\lambda t)^{\frac{x}{\lambda}}=\sum_{n=0}^\infty(x)_{n,\lambda}\frac{t^{n}}{n!}, \ \ e_\lambda(t)=e_\lambda^1(t),\quad (0 \ne \lambda\in\mathbb{R}),
\end{split}
\end{equation}

where

\begin{equation}\label{eq08}
\begin{split}
(x)_{0,\lambda}=1, \ (x)_{n,\lambda}=x(x-\lambda)(x-2\lambda) \cdots(x-(n-1)\lambda),\ \ (n\geq1), \quad \ \text{(see \ \cite{17})}.
\end{split}
\end{equation}

Let $\log_\lambda t$ be the compositional inverse of $e_\lambda (t)$ such that $\log_\lambda(e_\lambda(t))=e_\lambda (\log_\lambda t)=t$.

Then we have
\begin{equation}\label{eq09}
\begin{split}
\log_\lambda (1+t)=\sum_{n=1}^\infty \frac{\lambda^{n-1}(1)_{n,\frac{1}{\lambda}}}{n!}t^n, \quad \ \text{(see \ \cite{12})}.
\end{split}
\end{equation}

In view of \eqref{eq01} and \eqref{eq02}, the degenerate Stirling numbers of the first kind, $S_{1,\lambda}(n,k)$, and the degenerate Stirling numbers of the second kind, $S_{2,\lambda}(n,k)$, are defined by
\begin{equation}\label{eq10}
\begin{split}
(x)_n = \sum_{k=0}^n S_{1,\lambda}(n,k)(x)_{k,\lambda}, \ \quad (n\geq0),
\end{split}
\end{equation}

and

\begin{equation}\label{eq11}
\begin{split}
(x)_{n,\lambda}=\sum_{k=0}^n S_{2,\lambda}(n,k)(x)_{k}, \ \quad (n\geq0) \quad \  \text{(see \ \cite{12})}.
\end{split}
\end{equation}

Note that $\lim_{\lambda\rightarrow 0}\log_\lambda(1+t)=\log(1+t), \lim_{\lambda\rightarrow 0}e_\lambda^x(t)=e^{xt}$, $\lim_{\lambda\rightarrow 0}S_{1,\lambda}(n,k)=S_1 (n,k)$, and $\lim_{\lambda\rightarrow 0}S_{2,\lambda}(n,k)=S_2 (n,k)$, where $(n, k \geq 0)$.

From \eqref{eq10} and \eqref{eq11}, we note that
\begin{equation}\label{eq12}
\begin{split}
\frac{1}{k!}(e_\lambda (t)-1)^k = \sum_{n=k}^\infty S_{2,\lambda}(n,k) \frac{t^n}{n!}, \quad \ (k \geq 0),
\end{split}
\end{equation}

and

\begin{equation}\label{eq13}
\begin{split}
\frac{1}{k!}\big(\log_\lambda (1+t)\big)^k=\sum_{n=k}^\infty S_{1,\lambda}(n,k) \frac{t^n}{n!}, \quad \ \text{(see \ \cite{12})}.
\end{split}
\end{equation}

The degenerate Bell polynomials are defined by

\begin{equation}\label{eq14}
\begin{split}
e^{x(e_\lambda(t)-1)}=\sum_{n=0}^\infty \phi_{n,\lambda}(x) \frac{t^n}{n!}, \ \quad \text{(see \ \cite{12,15,16})}.
\end{split}
\end{equation}

Thus, by \eqref{eq12} and \eqref{eq14}, we get

\begin{equation}\label{eq15}
\begin{split}
\phi_{n,\lambda}(x)=\sum_{k=0}^n S_{2,\lambda}(n,k)x^k, \ \quad (n\geq 0), \ \quad \text{(see \ \cite{15,16,17})}.
\end{split}
\end{equation}

A random variable $X$ is a real valued function defined on a sample space. If $X$ takes any value in a countable set, then $X$ is called a discrete random variable. For a discrete random variable $X$, the probability mass function $p(a)$ of $X$ is defined by

\begin{equation}\label{eq16}
\begin{split}
p(a)=P\{X=a\}, \quad \ \ \text{(see \ \cite{10,11,20})}.
\end{split}
\end{equation}

A random variable $X$ taking on one of the values $0,1,2,\cdots$ is said to be the Poisson random variable with parameter $\alpha (>0)$, which is denoted by $X\sim Poi(\alpha)$, if the probability mass function of $X$ is given by
\begin{equation}\label{eq17}
\begin{split}
p(i)=P\{X=i\}=e^{-\alpha}\frac{\alpha^i}{i!}, \ i=0,1,2,\cdots , \ \quad \ \text{(see \ \cite{10,11,20})}.
\end{split}
\end{equation}

For $n \geq 1$, the quantity $E[X^n]$ of the Poisson random variable X with parameter $\alpha(>0)$, which is called the $n$th moment of $X$, is given by

\begin{equation}\label{eq18}
\begin{split}
E[X^n]=\sum_{i=0}^\infty i^np(i)=e^{-\alpha} \sum_{i=0}^\infty \frac{i^n}{i!}\alpha^i=\phi_n(\alpha), \ \quad \ \text{(see \ \cite{10,11,20})}.
\end{split}
\end{equation}

As is well known, the Charlier polynomials $C_n(x;\alpha)$ are defined by

\begin{equation}\label{eq19}
\begin{split}
e^{-\alpha t}(1+t)^x = \sum_{n=0}^\infty C_n(x;\alpha) \frac{t^n}{n!}, \ \ \quad \text{(see \ \cite{19})},
\end{split}
\end{equation}
where $x,t,\alpha \in \mathbb{R}$.

Thus, by \eqref{eq19}, we get
\begin{equation}\label{eq19-1}
\begin{split}
C_n(x;\alpha) = \sum_{l=0}^n \bigg(\sum_{k=l}^n \binom{n}{k}(-1)^{n-k}\alpha^{n-k} S_1(k,l)\bigg)x^l.
\end{split}
\end{equation}

A finite lattice $L$ is geometric if it is a finite semimodular lattice which is also atomic. Dowling constructed an important finite geometric lattice $Q_n(G)$ out of a finite set of $n$ elements and a finite group $G$ of order $m$, called Dowling lattice of rank $n$ over a finite group of order $m$. If $L$ is the Dowling lattice $Q_n(G)$ of rank $n$ over a finite group $G$ of order $m$, then the Whitney numbers of the first kind $V_{Q_n(G)}(n,k)$ and the Whitney numbers of the second kind $W_{Q_{n}(n,k)}$ are respectively denoted by $V_m(n,k)$ and $W_m(n,k)$. The Whitney numbers $V_m(n,k)$ and $W_m(n,k)$ satisfy the following Stirling number-like relations:

\begin{equation}\label{eq20}
\begin{split}
(mx+1)^n = \sum_{k=0}^n W_m(n,k)m^k(x)_k,
\end{split}
\end{equation}

\begin{equation}\label{eq21}
\begin{split}
m^n (x)_n = \sum_{k=0}^n V_m (n,k) (mx+1)^k,  \ \quad (n \geq 0), \ \quad \  \text{(see \ \cite{8, 9, 14})}.
\end{split}
\end{equation}

For $n \geq 0$, Dowling polynomials are given by

\begin{equation}\label{eq22}
\begin{split}
D_m(n,x)=\sum_{k=0}^n W_m(n,k)x^k, \ \ \ \text{(see \ \cite{8, 9, 13})}.
\end{split}
\end{equation}

Recently, Kim-Kim considered the degenerate Whitney numbers of the second kind defined by

\begin{equation}\label{eq23}
\begin{split}
(mx+1)_{n,\lambda}=\sum_{k=0}^n W_{m,\lambda}(n,k) m^k(x)_k, \ \quad (k \geq 0), \ \ \quad  \text{(see \ \cite{13})}.
\end{split}
\end{equation}

Thus, by \eqref{eq23}, we get

\begin{equation}\label{eq24}
\begin{split}
e_\lambda(t) \frac{1}{k!} \bigg(\frac{e_\lambda^m (t)-1}{m}\bigg)^k = \sum_{n=k}^\infty W_{m,\lambda}(n,k) \frac{t^n}{n!}.
\end{split}
\end{equation}

In \cite{14}, the degenerate Dowling polynomials are defined by

\begin{equation}\label{eq25}
\begin{split}
e_\lambda (t) e^{\frac{x}{m}(e_\lambda^m(t)-1)}=\sum_{n=0}^\infty D_{m,\lambda}(n,x) \frac{t^n}{n!}.
\end{split}
\end{equation}

By \eqref{eq24} and \eqref{eq25}, we get

\begin{equation}\label{eq26}
\begin{split}
D_{m,\lambda}(n,x)=\sum_{k=0}^n W_{m,\lambda}(n,k) x^k, \ \quad \text{(see \ \cite{13, 14})}.
\end{split}
\end{equation}

A further generalization of degenerate Whitney numbers of the second kind, Kim-Kim introduced the degenerate $r$-Whitney numbers of the second kind given by

\begin{equation}\label{eq27}
\begin{split}
(mx+r)_{n,\lambda}=\sum_{k=0}^n W_{m,\lambda}^{(r)}(n,k)m^k(x)_k, \ \quad (n,r \geq 0), \ \ \quad \text{(see \ \cite{13, 14})}.
\end{split}
\end{equation}

In view of \eqref{eq26}, they defined the degenerate $r$-Dowling polynomials given by

\begin{equation}\label{eq28}
\begin{split}
D_{m,\lambda}^{(r)}(n,x) = \sum_{k=0}^n W_{m,\lambda}^{(r)} (n,k)x^k, \ \quad (n \geq 0), \ \ \text{(see \ \cite{13, 14})}.
\end{split}
\end{equation}

From \eqref{eq28}, we can show that the generating function of the degenerate $r$-Dowling polynomials is given by

\begin{equation}\label{eq29}
\begin{split}
e_\lambda^r(t) e^{\frac{x}{m}(e_\lambda^m(t)-1)}=\sum_{n=0}^\infty D_{m,\lambda}^{(r)} (n,x) \frac{t^n}{n!}, \  \quad \text{(see \ \cite{13, 14})}.
\end{split}
\end{equation}

\section{Poisson degenerate central moments related to degenerate Dowling and degenerate $r$-Dowling polynomials}

Let $X$ be the Poisson random variable with mean $\frac{\alpha}{m}$. Then we consider the {\it{Poisson degenerate central moments}} given by $E[(mX+1)_{n,\lambda}], \  (n \geq 0)$. We observe that

\begin{equation}\label{eq30}
\begin{split}
E[e_\lambda^{mX+1}(t)]&=\sum_{k=0}^{\infty} e_\lambda^{mk+1}(t) p(k) \\
&=e_\lambda(t) e^{-\frac{\alpha}{m}} \sum_{k=0}^\infty e_\lambda^{mk}(t)\frac{1}{k!} \bigg(\frac{\alpha}{m}\bigg)^k \\
&=e_\lambda(t)e^{-\frac{\alpha}{m}} e^{\frac{\alpha}{m}e_\lambda^m(t)}=e_\lambda (t) e^{\frac{\alpha}{m}(e_\lambda^m(t)-1)} \\
&=\sum_{n=0}^\infty D_{m,\lambda}(n,\alpha)\frac{t^n}{n!}.
\end{split}
\end{equation}

On the other hand, by \eqref{eq07}, we get
\begin{equation}\label{eq31}
E[e_\lambda^{mX+1}(t)] =\sum_{n=0}^\infty E[(mX+1)_{n,\lambda}] \frac{t^n}{n!}.
\end{equation}

Therefore, by \eqref{eq30} and \eqref{eq31}, we obtain the following theorem.

\begin{theorem}
Let $X\sim Poi\big(\frac{\alpha}{m}\big)$. Then we have
\begin{equation*}
\begin{split}
E[(mX+1)_{n,\lambda}]=D_{m,\lambda}(n,\alpha), \ \quad (n\geq0).
\end{split}
\end{equation*}
\end{theorem}

\noindent Note that
\begin{equation*}
\begin{split}
E[(mX+1)^n]&=\lim_{\lambda\rightarrow 0}E[(mX+1)_{n,\lambda}]=\lim_{\lambda\rightarrow 0}D_{m,\lambda}(n,\alpha) =D_m(n,\alpha), \quad (n\geq0).
\end{split}
\end{equation*}

\medskip

Let $X\sim Poi\big(\frac{\alpha}{m}\big)$. It is not difficult to show that
\begin{equation}\label{eq32}
\begin{split}
(x+y)_{n,\lambda}=\sum_{k=0}^n\binom{n}{k}(x)_{k,\lambda}(y)_{n-k,\lambda},\quad \ (n\geq0).
\end{split}
\end{equation}

By \eqref{eq32}, we get
\begin{equation}\label{eq33}
\begin{split}
E[(mX+1)_{n,\lambda}]=\sum_{k=0}^n\binom{n}{k}(1)_{n-k,\lambda}E[(mX)_{k,\lambda}].
\end{split}
\end{equation}

From \eqref{eq17}, we have
\begin{equation}\label{eq34}
\begin{split}
\sum_{n=0}^\infty E[(mX)_{n,\lambda}]\frac{t^n}{n!}=E[e^{mX}_\lambda(t)]&=e^{-\frac{\alpha}{m}}\sum_{n=0}^\infty e_\lambda^{mn}(t)\frac{(\frac{\alpha}{m})^n}{n!}\\
&=e^{\frac{\alpha}{m}(e_\lambda^m(t)-1)}=e^{\frac{\alpha}{m}(e_{\frac{\lambda}{m}}(mt)-1)}=\sum_{n=0}^\infty\phi_{n,\frac{\lambda}{m}}\bigg(\frac{\alpha}{m}\bigg)m^n\frac{t^n}{n!}.
\end{split}
\end{equation}

Comparing the coefficients on both sides of \eqref{eq34}, we have
\begin{equation}\label{eq35}
\begin{split}
E[(mX)_{n,\lambda}]=\phi_{n,\frac{\lambda}{m}}\bigg(\frac{\alpha}{m}\bigg)m^n,\quad (n\geq0).
\end{split}
\end{equation}

By \eqref{eq15}, we get
\begin{equation}\label{eq36}
\begin{split}
\phi_{n,\frac{\lambda}{m}}\bigg(\frac{\alpha}{m}\bigg)=\sum_{k=0}^nS_{2,\frac{\lambda}{m}}(n,k)\bigg(\frac{\alpha}{m}\bigg)^k.
\end{split}
\end{equation}

Therefore, by \eqref{eq33}, \eqref{eq35} and \eqref{eq36}, we get
\begin{equation}\label{eq37}
\begin{split}
E[(mX+1)_{n,\lambda}]&=\sum_{k=0}^n\binom{n}{k}(1)_{n-k,\lambda}m^k\phi_{k,\frac{\lambda}{m}}\bigg(\frac{\alpha}{m}\bigg)\\
&=\sum_{k=0}^n\binom{n}{k}(1)_{n-k,\lambda}m^k\sum_{j=0}^k S_{2,\frac{\lambda}{m}}(k,j)\bigg(\frac{\alpha}{m}\bigg)^j\\
&=\sum_{j=0}^n \alpha^j \sum_{k=j}^n\binom{n}{k}(1)_{n-k,\lambda}m^{k-j}S_{2,\frac{\lambda}{m}}(k,j).
\end{split}
\end{equation}

Therefore, by \eqref{eq37}, we obtain the following theorem.

\begin{theorem}

For $n\geq0$, let $X\sim Poi(\frac{\alpha}{m})$. Then we have
\begin{equation*}
\begin{split}
E[(mX+1)_{n,\lambda}]&=\sum_{k=0}^n\binom{n}{k}(1)_{n-k,\lambda}m^k\phi_{k,\frac{\lambda}{m}}\bigg(\frac{\alpha}{m}\bigg)\\
&=\sum_{j=0}^n\alpha^j\sum_{k=j}^n\binom{n}{k}(1)_{n-k,\lambda}m^{k-j}S_{2,\frac{\lambda}{m}}(k,j).
\end{split}
\end{equation*}

\end{theorem}

\medskip

When $m=1$, we have
\begin{equation*}
\begin{split}
\sum_{n=0}^\infty D_{1,\lambda}(n,\alpha)\frac{t^n}{n!}&=\sum_{n=0}^\infty E[(X+1)_{n,\lambda}]\frac{t^n}{n!}=E[e_\lambda^{X+1}(t)]\\
&=e_\lambda(t)e^{-\alpha}\sum_{n=0}^\infty\frac{e_\lambda^n(t)}{n!}\alpha^n=e_\lambda(t)e^{\alpha(e_\lambda(t)-1)}\\
&=e^{\alpha(e_\lambda(t)-1)}+\frac{d}{d\alpha}e^{\alpha(e_\lambda(t)-1)}=\sum_{n=0}^\infty\Big(\phi_{n,\lambda}(\alpha)+\frac{d}{d\alpha}\phi_{n,\lambda}(\alpha)\Big)\frac{t^n}{n!}.
\end{split}
\end{equation*}

Thus, we have
\begin{equation}\label{eq38}
\begin{split}
\phi_{n,\lambda}(\alpha)+\frac{d}{d\alpha}\phi_{n,\lambda}(\alpha)=D_{1,\lambda}(n,\alpha)=E[(X+1)_{n,\lambda}],\quad (n\geq0).
\end{split}
\end{equation}


On the other hand, by \eqref{eq15}, we have
\begin{equation}\label{eq39}
\begin{split}
\frac{d}{d\alpha}\phi_{n,\lambda}(\alpha)=\frac{d}{d\alpha}\sum_{k=0}^n S_{2,\lambda}(n, \ k)\alpha^k&=\sum_{k=1}^n k S_{2,\lambda}(n, \ k)\alpha^{k-1}\\
&=\sum_{k=0}^{n-1}(k+1)S_{2,\lambda}(n, \ k+1)\alpha^k.
\end{split}
\end{equation}

Therefore, by \eqref{eq39}, we obtain the following corollary.

\begin{corollary}

For $n\geq0$, let $X\sim Poi(\alpha)$. Then we have
\begin{equation*}
\begin{split}
E[(X+1)_{n,\lambda}]=D_{1,\lambda}(n,\ \alpha)=\sum_{k=0}^{n}\Big((k+1)S_{2,\lambda}(n, \ k+1)+S_{2,\lambda}(n,k)\Big)\alpha^k.
\end{split}
\end{equation*}

\end{corollary}

\medskip

For $r\geq0$, let $X$ be the Poisson random variable with mean $\frac{\alpha}{m}$. Then we have
\begin{equation}\label{eq40}
\begin{split}
E[e_\lambda^{mX+r}(t)]&=\sum_{k=0}^\infty e_\lambda^{mk+r}(t)p(k)\\
&=e_\lambda^r(t)e^{-\frac{\alpha}{m}}\sum_{k=0}^\infty\frac{e_\lambda^{mk}(t)}{k!}\bigg(\frac{\alpha}{m}\bigg)^k\\
&=e_\lambda^r(t)e^{-\frac{\alpha}{m}}e^{\frac{\alpha}{m}e_\lambda^m(t)}=e_\lambda^r(t)e^{\frac{\alpha}{m}(e_\lambda^m(t)-1)}=\sum_{n=0}^\infty D_{m,\lambda}^{(r)}(n, \ \alpha)\frac{t^n}{n!}.
\end{split}
\end{equation}


The left hand side of \eqref{eq40} is given by
\begin{equation}\label{eq41}
\begin{split}
E[e_\lambda^{mX+r}(t)]=\sum_{n=0}^\infty E[(mX+r)_{n,\lambda}]\frac{t^n}{n!}.
\end{split}
\end{equation}

Therefore, by \eqref{eq40} and \eqref{eq41}, we obtain the following theorem.
\begin{theorem}

For $m,\ r\geq0$, let $X\sim Poi(\frac{\alpha}{m})$. Then we have \begin{equation*}
\begin{split}
E[(mX+r)_{n,\lambda}]=D_{m,\lambda}^{(r)}(n,\alpha),\quad (n\geq0).
\end{split}
\end{equation*}

\end{theorem}

\noindent Note that
\begin{equation*}
\begin{split}
D_m^{(r)}(n,\alpha)=\lim_{\lambda\rightarrow0}E[(mX+r)_{n,\lambda}]=E[(mX+r)^n].
\end{split}
\end{equation*}

By \eqref{eq32}, we get
\begin{equation}\label{eq42}
\begin{split}
E[(mX+r)_{n,\lambda}]&=\sum_{k=0}^n\binom{n}{k}(r)_{n-k,\lambda}E[(mX)_{k,\lambda}]\\
&=\sum_{k=0}^n\binom{n}{k}(r)_{n-k,\lambda}m^k\phi_{k,\frac{\lambda}{m}}\bigg(\frac{\alpha}{m}\bigg). \end{split}
\end{equation}

Therefore, by \eqref{eq42}, we obtain the following corollary.
\begin{corollary}

For $n, \ r\geq0$, let $X\sim Poi(\frac{\alpha}{m})$. Then we have
\begin{equation*}
\begin{split}
D_{m,\lambda}^{(r)}(n,\alpha)=E[(mX+r)_{n,\lambda}]=\sum_{k=0}^n\binom{n}{k}(r)_{n-k,\lambda}m^k\phi_{k,\frac{\lambda}{m}}\bigg(\frac{\alpha}{m}\bigg).
\end{split}
\end{equation*}

\end{corollary}

When $r=1$, we have
\begin{equation*}
\begin{split}
E[(mX+1)_{n,\lambda}]=D_{m,\lambda}^{(1)}(n,\alpha)=D_{m,\lambda}(n,\alpha).
\end{split}
\end{equation*}

\medskip

From \eqref{eq27}, we have
\begin{equation}\label{eq43}
\begin{split}
D_{m,\lambda}^{(r)}(n,\alpha)=E[(mX+r)_{n,\lambda}]=\sum_{k=0}^nW_{m,\lambda}^{(r)}(n,k)m^kE[(X)_k].
\end{split}
\end{equation}

By \eqref{eq03}, we get
\begin{equation}\label{eq44}
\begin{split}
\sum_{n=0}^\infty E[(X)_n]\frac{t^n}{n!}&=E[(1+t)^X]=\sum_{k=0}^\infty E[X^k]\frac{1}{k!}(\log(1+t))^k\\
&=\sum_{k=0}^\infty E[X^k]\sum_{n=k}^\infty S_1(n,k)\frac{t^n}{n!}=\sum_{n=0}^\infty\bigg(\sum_{k=0}^n S_1(n,k)E[X^k]\bigg)\frac{t^n}{n!}.
\end{split}
\end{equation}

Since $X\sim Poi(\frac{\alpha}{m})$, from \eqref{eq17}, we have
\begin{equation}\label{eq45}
\begin{split}
\sum_{k=0}^\infty E[X^k]\frac{t^k}{k!}=E[e^{Xt}]&=e^{-\frac{\alpha}{m}}\sum_{k=0}^\infty e^{kt}\frac{(\frac{\alpha}{m})^k}{k!}\\
&=e^{\frac{\alpha}{m}(e^t-1)}=\sum_{k=0}^\infty \phi_k\bigg(\frac{\alpha}{m}\bigg)\frac{t^k}{k!}. \end{split}
\end{equation}

Thus, by \eqref{eq45}, we get
\begin{equation}\label{eq46}
\begin{split}
E[X^k]=\phi_k\bigg(\frac{\alpha}{m}\bigg),\quad (k\geq0).
\end{split}
\end{equation}


From \eqref{eq44} and \eqref{eq46}, we have
\begin{equation}\label{eq47}
\begin{split}
E[(X)_n]=\sum_{k=0}^nS_1(n,k)E[X^k]=\sum_{k=0}^nS_1(n,k)\phi_k\bigg(\frac{\alpha}{m}\bigg).
\end{split}
\end{equation}

By \eqref{eq43} and \eqref{eq47}, we get
\begin{equation}\label{eq48}
\begin{split}
D_{m,\lambda}^{(r)}(n,\alpha)&=E[(mX+r)_{n,\lambda}]=\sum_{k=0}^n W_{m,\lambda}^{(r)}(n,k)m^kE[(X)_k]\\
&=\sum_{k=0}^nW_{m,\lambda}^{(r)}(n,k)m^k\sum_{j=0}^kS_1(k,j)\phi_j\bigg(\frac{\alpha}{m}\bigg)\\
&=\sum_{j=0}^n\bigg(\sum_{k=j}^n W_{m,\lambda}^{(r)}(n,k)m^k S_1(k,j)\bigg)\phi_j\bigg(\frac{\alpha}{m}\bigg).
\end{split}
\end{equation}

Therefore, by \eqref{eq48}, we obtain the following theorem.

\begin{theorem}

For $n\geq0$, let $X\sim Poi(\frac{\alpha}{m})$. Then we have
\begin{equation*}
\begin{split}
D_{m,\lambda}^{(r)}(n,\alpha)&=E[(mX+r)_{n,\lambda}]\\
&=\sum_{j=0}^n\bigg(\sum_{k=j}^nW_{m,\lambda}^{(r)}(n,k)m^k S_1(k,j)\bigg)\phi_j\bigg(\frac{\alpha}{m}\bigg).
\end{split}
\end{equation*}

\end{theorem}


\medskip

By Theorem 2, we get
\begin{equation}\label{eq49}
\begin{split}
D_{m,\lambda}^{(r)}(n,\alpha)&=E[(mX+r)_{n,\lambda}]=\sum_{k=0}^n\binom{n}{k}(r-1)_{n-k,\lambda}E[(mX+1)_{k,\lambda}]\\
&=\sum_{k=0}^n\binom{n}{k}(r-1)_{n-k,\lambda}\sum_{j=0}^k\alpha^j\sum_{l=j}^k\binom{k}{l}(1)_{k-l,\lambda}m^{l-j}S_{2,\frac{\lambda}{m}}(l,j)\\
&=\sum_{j=0}^n\alpha^j\bigg\{\sum_{k=j}^n\sum_{l=j}^k\binom{n}{k}\binom{k}{l}(r-1)_{n-k,\lambda}(1)_{k-l,\lambda}m^{l-j}S_{2,\frac{\lambda}{m}}(l,j)\bigg\}.
\end{split}
\end{equation}

On the other hand, by \eqref{eq28}, we get
\begin{equation}\label{eq50}
\begin{split}
D_{m,\lambda}^{(r)}(n,\ \alpha)=\sum_{j=0}^n\alpha^j W_{m,\lambda}^{(r)}(n,j).
\end{split}
\end{equation}

Therefore, by \eqref{eq49} and \eqref{eq50}, we obtain the following theorem.
\begin{theorem}

For $n,\ j\geq0$, we have
\begin{equation*}
\begin{split}
W_{m,\lambda}^{(r)}(n,j)=\sum_{k=j}^n\sum_{l=j}^k\binom{n}{k}\binom{k}{l}(r-1)_{n-k,\lambda}(1)_{k-l,\lambda}m^{l-j}S_{2,\frac{\lambda}{m}}(l,j).
\end{split}
\end{equation*}

\end{theorem}

\medskip

From \eqref{eq50} and \eqref{eq42}, we note that
\begin{equation}\label{eq51}
\begin{split}
\sum_{j=0}^n\alpha^j W_{m,\lambda}^{(r)}(n,j)&=D_{m,\lambda}^{(r)}(n,\alpha)=E[(mX+r)_{n,\lambda}]\\
&=\sum_{l=0}^n\binom{n}{l}(r)_{n-l,\lambda}m^l\phi_{l,\frac{\lambda}{m}}\bigg(\frac{\alpha}{m}\bigg)\\
&=\sum_{l=0}^n\binom{n}{l}(r)_{n-1,\lambda}m^l\sum_{j=0}^l\bigg(\frac{\alpha}{m}\bigg)^j S_{2,\frac{\lambda}{m}}(l,j)\\
&=\sum_{j=0}^n\alpha^j\bigg(\sum_{l=j}^n\binom{n}{l}(r)_{n-l,\lambda}m^{l-j}S_{2,\frac{\lambda}{m}}(l,j)\bigg).
\end{split}
\end{equation}


Therefore, by comparing the coefficients on both sides of \eqref{eq51}, we obtain the following theorem.

\begin{theorem}

For $n,j\geq0$, we have
\begin{equation*}
\begin{split}
W_{m,\lambda}^{(r)}(n,j)=\sum_{l=j}^n\binom{n}{l}(r)_{n-l,\lambda}m^{l-j}S_{2,\frac{\lambda}{m}}(l,j).
\end{split}
\end{equation*}

\end{theorem}

\medskip

We recall that Charlier polynomials $C_n(x;\alpha)$ are given by
\begin{equation}\label{eq52}
\begin{split}
e^{-\alpha t}(1+t)^x=\sum_{n=0}^\infty C_n(x\ ;\alpha)\frac{t^n}{n!},\quad \ (n\geq0 \ \ {\rm{and}} \ \ x,\alpha, t\in \mathbb{R}).
\end{split}
\end{equation}

Let us take $x=0$. Then we have
\begin{equation}\label{eq53}
\begin{split}
e^{-\alpha t}=\sum_{n=0}^\infty C_n(0\ ;\alpha)\frac{t^n}{n!}.
\end{split}
\end{equation}

Replacing $t$ by $1-e_\lambda^m(t)$, we get \begin{equation}\label{eq54}
\begin{split}
e^{\alpha(e_\lambda^m(t)-1)}=\sum_{k=0}^\infty C_k(0\ ;\alpha)(-1)^k\frac{1}{k!}(e_\lambda^m(t)-1)^{k}.
\end{split}
\end{equation}

Thus, by \eqref{eq14} and \eqref{eq54}, we get
\begin{equation}\label{eq55}
\begin{split}
\sum_{n=0}^\infty\phi_{n,\frac{\lambda}{m}}(\alpha)m^n\frac{t^n}{n!}&=e^{\alpha(e_{\frac{\lambda}{m}}(mt)-1)}=e^{\alpha(e_\lambda^m(t)-1)}\\
&=\sum_{k=0}^\infty C_k(0\ ;\alpha)(-1)^k\frac{1}{k!}(e_{\frac{\lambda}{m}}(mt)-1)^k\\&=\sum_{k=0}^\infty C_k(0\ ;\alpha)(-1)^k\sum_{n=k}^\infty S_{2,\frac{\lambda}{m}}(n,k)m^n\frac{t^n}{n!}\\
&=\sum_{n=0}^\infty\bigg(m^n\sum_{k=0}^nC_k(0\ ; \alpha)(-1)^k S_{2,\frac{\lambda}{m}}(n,k)\bigg)\frac{t^n}{n!}.
\end{split}
\end{equation}

Therefore, by \eqref{eq55}, we obtain the following theorem.
\begin{theorem}

For $n\geq0$, let $X\sim Poi(\alpha)$, Then we have
\begin{equation*}
\begin{split}
\phi_{n,\frac{\lambda}{m}}(\alpha)=E[(X)_{n,\frac{\lambda}{m}}]=\sum_{k=0}^nC_k(0\ ; \alpha)(-1)^k S_{2,\frac{\lambda}{m}}(n,k). \end{split}
\end{equation*}

\end{theorem}

\medskip

Let us take $x=1$ in \eqref{eq52}. Then we have
\begin{equation}\label{eq56}
\begin{split}
e^{-\alpha t}(1+t)=\sum_{n=0}^\infty C_n(1\ ;\alpha)\frac{t^n}{n!}.
\end{split}
\end{equation}

Replacing $t$ by $e_\lambda^m(t)-1$ and $\alpha$ by $-\frac{\alpha}{m}$, we get
\begin{equation}\label{eq57}
\begin{split}
\sum_{n=0}^\infty D_{m,\lambda}^{(m)}(n,\alpha)\frac{t^{n}}{n!}=e_\lambda^m(t)&e^{\frac{\alpha}{m}(e_\lambda^m(t)-1)}=\sum_{k=0}^\infty C_k\bigg(1\ ; -\frac{\alpha}{m}\bigg)\frac{(e_\lambda^m(t)-1)^k}{k!}\\
&=\sum_{k=0}^\infty C_k\bigg(1\ ; -\frac{\alpha}{m}\bigg)\frac{1}{k!}(e_{\frac{\lambda}{m}}(mt)-1)^k\\&=\sum_{n=0}^\infty\bigg(\sum_{k=0}^nC_k\bigg(1\ ; -\frac{\alpha}{m}\bigg)m^nS_{2,\frac{\lambda}{m}}(n,k)\bigg)\frac{t^n}{n!}.
\end{split}
\end{equation}

Comparing the coefficients on both sides of \eqref{eq57}, we have
\begin{equation}\label{eq58}
\begin{split}
D_{m,\lambda}^{(m)}(n,\alpha)=\sum_{k=0}^nC_k\bigg(1\ ;-\frac{\alpha}{m}\bigg)m^n S_{2,\frac{\lambda}{m}}(n,k),\quad (n\geq0).
\end{split}
\end{equation}

Note that
\begin{equation}\label{eq59}
\begin{split}
E[(mX+r)_{n,\lambda}]&=E[(mX+m+r-m)_{n,\lambda}]\\
&=\sum_{k=0}^n\binom{n}{k}(r-m)_{n-k,\lambda}E[(mX+m)_{k,\lambda}]\\
&=\sum_{k=0}^n\binom{n}{k}(r-m)_{n-k,\lambda}\sum_{j=0}^kC_j\bigg(1\ ;-\frac{\alpha}{m}\bigg)S_{2,\frac{\lambda}{m}}(k,j)m^k\\
&=\sum_{j=0}^nC_j\bigg(1\ ; -\frac{\alpha}{m}\bigg)\sum_{k=j}^n\binom{n}{k}(r-m)_{n-k,\lambda}S_{2,\frac{\lambda}{m}}(k,j)m^k.
\end{split}
\end{equation}

Therefore, by \eqref{eq58} and \eqref{eq59}, we obtain the following theorem.

\begin{theorem}
For $n \geq 0$, let $X\sim Poi (\frac{\alpha}{m})$. Then we have
\begin{equation*}
\begin{split}
E[(mX+m)_{n,\lambda}]=D_{m,\lambda}^{(m)}(n,\alpha) = \sum_{k=0}^n C_k(1;-\frac{\alpha}{m})m^n S_{2,\frac{\lambda}{m}} (n,k).
\end{split}
\end{equation*}
Furthermore, we have
\begin{equation*}
\begin{split}
D_{m,\lambda}^{(r)}(n,\alpha) &= E[(mX+r)_{n,\lambda}] \\
&= \sum_{j=0}^n C_j \bigg(1;-\frac{\alpha}{m}\bigg)\sum_{k=j}^n \binom{n}{k}(r-m)_{n-k,\lambda} S_{2,\frac{\lambda}{m}}(k,j)m^k.
\end{split}
\end{equation*}
\end{theorem}

\medskip

Replacing $t$ by $(e_\lambda^m(t)-1)$, $\alpha$ by $-\frac{\alpha}{m}$, and $x$ by $\frac{r}{m}$ in \eqref{eq52}, we have

\begin{equation}\label{eq60}
\begin{split}
e^{\frac{\alpha}{m}(e_\lambda^m(t)-1)}e_\lambda ^r(t) &= \sum_{k=0}^\infty C_k \bigg(\frac{r}{m} ; -\frac{\alpha}{m}\bigg) \frac{1}{k!}(e_\lambda^m(t)-1)^k \\
&=\sum_{k=0}^\infty C_k\bigg(\frac{r}{m} ; -\frac{\alpha}{m}\bigg) \frac{1}{k!}(e_{\frac{\lambda}{m}}(mt)-1)^k\\& =\sum_{k=0}^\infty C_k\bigg(\frac{r}{m} ; -\frac{\alpha}{m}\bigg)\sum_{n=k}^\infty S_{2,\frac{\lambda}{m}}(n,k) m^n \frac{t}{n!} \\
&=\sum_{n=0}^\infty \bigg(m^n \sum_{k=0}^n C_k\bigg(\frac{r}{m} ; -\frac{\alpha}{m}\bigg)S_{2,\frac{\lambda}{m}}(n,k)\bigg)\frac{t^n}{n!}.
\end{split}
\end{equation}

On the other hand, by \eqref{eq29}, we get
\begin{equation}\label{eq61}
\begin{split}
e_\lambda^r(t)e^{\frac{\alpha}{m}(e_\lambda^m(t)-1)}=\sum_{n=0}^\infty D_{m,\lambda}^{(r)}(n,\alpha)\frac{t^n}{n!}.
\end{split}
\end{equation}

Therefore, by \eqref{eq60} and \eqref{eq61}, we obtain the following theorem.

\begin{theorem}

For $n\geq0$, let $X\sim Poi(\frac{\alpha}{m})$. Then we have
\begin{equation*}
\begin{split}
D_{m,\lambda}^{(r)}(n,\alpha)=E[(mX+r)_{n,\lambda}]=m^n\sum_{k=0}^nC_k\bigg(\frac{r}{m}\ ;-\frac{\alpha}{m} \bigg) S_{2,\frac{\lambda}{m}}(n,k).
\end{split}
\end{equation*}

\end{theorem}

\medskip

\section{Conclusion}
In recent years, studying various degenerate versions of many special polynomials and numbers received regained interests of some mathematicians and many interesting results were discovered.
Degenerate Dowling and degenerate $r$-Dowling polynomials were introduced earlier as degenerate versions and further generalizations of Dowling and $r$-Dowling polynomials. \par
Assume that $X$ is the Poisson random variable with mean $\frac{\alpha}{m}$. We showed that the Poisson degenerate central moment $E[(mX+1)_{n,\lambda}]$ is equal to $D_{m,\lambda}(n,\alpha)$ and to an expression involving  the degenerate Bell polynomials, respectively in Theorem 1 and Theorem 2. We deduced that $E[(mX+r)_{n,\lambda}]$ is equal to $D_{m,\lambda}^{(r)}(n,\alpha)$ in Theorem 4. We expressed the same in terms of the degenerate Bell polynomials in Corollary 5 and of the degenerate $r$-Whitney numbers of the second kind and the Bell polynomials in Theorem 6. Furthermore, it is represented by the Charlier polynomials and the degenerate Stirling numbers of the second kind in Theorems 10 and 11. \par
As one of our future projects, we would like to continue to study degenerate versions of certain special polynomials and numbers and their applications to physics, science and engineering as well as mathematics.

\bigskip

\noindent{\bf{Acknowledgments}} \\
The authors thank Jangjeon Institute for Mathematical Sciences for the support of this research.

\vspace{0.1in}

\noindent{\bf {Availability of data and material}} \\
Not applicable.

\vspace{0.1in}

\noindent{\bf{Funding}} \\
This work was supported by the Basic Science Research Program, the National
               Research Foundation of Korea,
               (NRF-2021R1F1A1050151).
\vspace{0.1in}

\noindent{\bf{Ethics approval and consent to participate}} \\
The authors declare that there is no ethical problem in the production of this paper.

\vspace{0.1in}

\noindent{\bf {Competing interests}} \\
The authors declare no conflict of interest.

\vspace{0.1in}

\noindent{\bf{Consent for publication}} \\
The authors want to publish this paper in this journal.

\end{document}